# The data-driven Schrödinger bridge

Michele Pavon[*]     Esteban G. Tabak[†]     Giulio Trigila [‡]

June 3, 2018

## Abstract

Erwin Schrödinger posed, and to a large extent solved in 1931/32 the problem of finding the most likely random evolution between two continuous probability distributions. This article considers this problem in the case when only samples of the two distributions are available. A novel iterative procedure is proposed, inspired by Fortet-Sinkhorn type algorithms. Since only samples of the marginals are available, the new approach features constrained maximum likelihood estimation in place of the nonlinear boundary couplings, and importance sampling to propagate the functions $\varphi$ and $\hat{\varphi}$ solving the Schrödinger system. This method is well-suited to high-dimensional settings, where introducing grids leads to numerically unfeasible or unreliable methods. The methodology is illustrated in two applications: entropic interpolation of two-dimensional Gaussian mixtures, and the estimation of integrals through a variation of importance sampling.

## 1   Introduction

This article proposes a methodology for solving the following problem: given $m$ and $n$ independent samples $\{x_i\}$ and $\{y_j\}$ from two distributions with probability densities $\rho_0(x)$ and $\rho_1(y)$ respectively, and a prior probability $p(t_1, x, t_2, y)$ that a "particle" at position $x$ at time $t_1$ will end up at position $y$ at time $t_2$, find the most likely intermediate evolution $\rho(z,t)$, $t \in [0,1]$ satisfying $\rho(x,0) = \rho_0(x)$ and $\rho(y,1) = \rho_1(y)$. This is a data-driven version of the Schrödinger bridge problem, which we describe below. In addition to the evolving density $\rho(z,t)$, the solution provides the posterior transition density $p^*(t_1, x, t_2, y)$ most consistent with the observed initial and final distributions, useful for model improvement.

---

[*]M. Pavon is with the Dipartimento di Matematica "Tullio Levi-Civita", Università di Padova, via Trieste 63, 35121 Padova, Italy; pavon@math.unipd.it

[†]E.G. Tabak is with the Courant Institute of Mathematical Sciences, 251 Mercer St., New York University, New York, NY 10012; tabak@cims.nyu.edu

[‡]G. Trigila is with the Department of Mathematics, Baruch College, 55 Lexington Ave, New York, NY10010; Giulio.Trigila@baruch.cuny.edu

## 1.1 Motivation, examples and extensions

Many problems of practical and theoretical interest can be directly formulated as data-driven Schrödinger bridges. Consider the following two examples, arising in climate studies and evolutionary biology:

1. With the current knowledge of oceanic or atmospheric flows described in terms of a velocity field $v(x,t)$ and a diffusion operator $D$, the corresponding Fokker-Planck evolution equation yields the prior $p(t_1, x, t_2, y)$ for the trajectories of tracers. If at any point in time a cloud of particles is released into the fluid (a volcanic eruption, a designed experiment) or its current concentration $\rho_0$ is sampled, and at some other time its distribution $\rho_1$ is sampled again, the data-driven bridge problem provides an estimate for the most likely intermediate evolution $\rho(z,t)$ of the tracer cloud and to an improved model for the currents $v$.

2. Given the distribution of traits (genomic or phenomic) for a species at two points in time, and a stochastic model for their evolution, the problem asks for the most likely intermediate evolutionary stages, and provides as additional output an improved stochastic evolutionary model.

In other problems, it is not an intermediate evolution that one is after, but the probabilistic matching $\pi(x, y)$ between two distributions $\rho_0(x)$ and $\rho_1(y)$ under a prior matching model $p(y|x)$. In this case, both the problem and the methodology proposed for solving it proposed extend without changes to situations where the variables $x$ and $y$ do not have the same dimensions, arising frequently in practice. For instance, in applications to the employment market, there is no reason for the number of variables characterizing employers and employees to be the same.

In a third type of scenarios, there is only one data-given distribution $\rho_1(y)$; the other distribution $\rho_0(x)$ and the prior $p(t_1, x, t_2, y)$ are introduced for convenience by the modeler, so as to perform $\rho_1(y)$-related tasks. As an example, in Section 5 we apply the Schrödinger bridge to develop a variation of importance sampling where the distribution over which expected values of a function are sought is known only through samples.

In other applications, one has only $\rho_1(y)$ and the prior $p(t_1, x, t_2, y)$, and would like to determine $\rho_t(z)$ for $t < 1$. Two prototypal examples are inverse problems, such as describing the most likely previous temperature distribution of a system given its current one, and large deviation problems: if the stochastic process described by $p$ has a statistically steady state $\rho_{eq}$, what are the most likely paths that will lead to a $\rho_1$ different from $\rho_{eq}$, such as the one corresponding to a strong storm or a draught in applications to weather and climate.

## 1.2   The methodology

The solution to the Schrödinger bridge problem can be factorized in the form (see (18) in Section 2 below)

$$\rho_t(x) = \varphi(t, x)\hat{\varphi}(t, x),$$

where $\rho_t(x)$ represents the distribution at time $t$, and $\varphi$ and $\hat{\varphi}$ evolve from $t = 1$ and $t = 0$ respectively, following the prior:

$$\hat{\varphi}(t, y) = \int p(0, x, t, y)\hat{\varphi}(0, x)dx,$$

$$\varphi(t, x) = \int p(t, x, 1, y)\varphi(1, y)dy.$$

One can therefore, starting from an arbitrary $\hat{\varphi}(0, x)$, propagate it into the corresponding $\hat{\varphi}(1, y)$, and write

$$\varphi(1, y) = \frac{\rho_1(y)}{\hat{\varphi}(1, y)}.$$

Then, evolving $\varphi(1, y)$ back into the corresponding $\varphi(0, x)$, we write

$$\hat{\varphi}(0, x) = \frac{\rho_0(x)}{\varphi(0, x)},$$

and repeat. This idea underlies iteration schemes that, under suitable assumptions, converge to the solution of the Schrödinger bridge problem [33, 13].

Yet this procedure assumes that the initial and final distributions $\rho_0$ and $\rho_1$, as well as the transition probability $p$, are known explicitly, and that the integrals propagating $\varphi$ and $\hat{\varphi}$ between $t = 0$ and $t = 1$ can be evaluated in closed form. By contrast, in applications $\rho_0$ and $\rho_1$ are typically only known through samples. In addition, it is often the case that the transition probability $p$ can be sampled through the integration of a stochastic differential equation, but not evaluated, which would require solving the corresponding Fokker-Plank equation. Moreover, even if $\rho_0$, $\rho_1$ and $p$ are known, one still needs to estimate the integrals propagating $\varphi$ and $\hat{\varphi}$ numerically.

The methodology developed in this article mimics the iterative procedure above, but replacing each step by a sample-based equivalent. Thus the statements that

$$\varphi(0,x)\hat{\varphi}(0,x) = \rho_0(x) \quad \text{and} \quad \varphi(1,y)\hat{\varphi}(1,y) = \rho_1(y)$$

are interpreted as density estimations and implemented via maximum likelihood, and the propagators for $\varphi$ and $\hat{\varphi}$ are estimated via importance sampling. Both tasks involve elements unique to the Schrödinger bridge problem, described in Section 4.

## 1.3 Prior work

Schrödinger's statistical mechanical thought experiment (large deviations problem) was motivated by analogies with quantum mechanics. On the other hand, since Boltzmann's fundamental work [8], and then through Sanov's theorem [51], we know that finding the most likely *Zustandverteilung* (macrostate) is equivalent to solving a maximum entropy problem. This connection provides a second important motivation for Schrödinger bridges, as an inference methodology that prescribes a posterior distribution making the fewest number of assumptions beyond the available information. This approach has been developed over the years, thanks in particular to the work of Jaynes, Burg, Dempster and Csiszár [39, 40, 9, 10, 27, 20, 21, 22]. A more recent third motivation for studying Schrödinger bridges is that they can be viewed as regularization of the Optimal Mass Transport (OMT) problem [44, 45, 46, 41, 42, 12] which mitigates its computational challenges [2, 3, 50]. A large number of papers have since appeared on computational regularized OMT (Sinkhorn-type algorithms), see e.g. [23, 4, 15, 13, 18, 43, 1, 19]. While most of the classical work concentrates on the continuous problem, see e.g. the bibliography in [42] and Section 2 below, these papers concern the discrete Schödinger bridge problem [48, 34]. Hardly any attention, however, has been given to the case when only samples of continuous marginals are available (one exception is [28] which deals with using regularized optimal transport for hard and soft clustering). One might think that the latter case may be readily treated by discretizing the spatial variables through grids. As we argue in the beginning of Section 4, such an approach is often numerically unfeasible and/or not reliable. Thus, in this paper we provide what appears to be the first numerically viable approach to the data-driven continuous Schrödinger Bridge problem.

As discussed at the end of Subsection 2.5, this approach permits finding a map from $\rho_0$ to $\rho_1$, relating this work to [58] and [56] developed in the context of optimal transport.

### 1.4 Organization of the article

The paper is organized as follows. In Section 2, we provide an introduction to Schrödinger Bridges. We include a concise description of Schrödinger's original motivation, and elements of the connection between the large deviation problem and a path space maximum entropy problem, and with Optimal Transport. We also sketch derivations of the Schrödinger system and of the stochastic control and fluid dynamic formulations, focusing on the case when the prior transition density is the heat kernel.

In Section 3, we outline Fortet's iterative algorithm, dating back to 1940, which represents a sort of guideline for the numerical methods we develop in the rest of the paper. Section 4 features the novel methodology to attack the data-driven bridge problem, motivated by numerical, statistical and optimization considerations. First, the so-called half-bridge problem is treated, and then the full bridge, leading to the algorithm of Subsection 4.3. In Section 5, we illustrate the methodology in two relevant applications: the entropic interpolation between two Gaussian mixtures on $\mathbb{R}^2$ and a new application of Schrödinger Bridges to a variation of Importance Sampling. Finally, in Section 6 we summarize the results and propose future avenues of research.

## 2 Background on Schrödinger Bridges

### 2.1 Schrödinger's hot gas experiment

In 1931/32, Erwin Schrödinger proposed the following *Gedankenexperiment* [52, 53]. Consider the evolution of a cloud of $N$ independent Brownian particles in $\mathbb{R}^n$. This cloud of particles has been observed having at the initial time $t = 0$ an empirical distribution equal to $\rho_0(x)dx$. At time $t = 1$, an empirical distribution equal to $\rho_1(x)dx$ is observed which considerably differs from what it should be according to the law of large numbers ($N$ is large, typically of the order of Avogadro's number), namely

$$\rho_1(y) \neq \int_{\mathbb{R}^3} p(0, x, 1, y)\rho_0(x)dx,$$

where

$$p(s, y, t, x) = [2\pi(t-s)]^{-\frac{n}{2}} \exp\left[-\frac{|x-y|^2}{2(t-s)}\right], \quad s < t \qquad (1)$$

is the transition density of the Wiener process. It is apparent that the particles have been transported in an unlikely way. But of the many unlikely ways in which this could have happened, which one is the most likely? In modern probabilistic terms, this is a problem of *large deviations of the empirical distribution* as observed by Föllmer [32].

## 2.2 Large deviations and maximum entropy formulation

The area of large deviations is concerned with the probabilities of very rare events. Thanks to Sanov's theorem [51], Schrödinger's problem can be turned into a maximum entropy problem for distributions on trajectories. Let $\Omega = C([0, 1]; \mathbb{R}^n)$ be the space of $\mathbb{R}^n$ valued continuous functions and let $X^1, X^2, \ldots$ be i.i.d. Brownian evolutions on $[0, 1]$ with values in $\mathbb{R}^n$ ($X_i$ is distributed according to the Wiener measure $W$ on $C([0, 1]; \mathbb{R}^n)$). The *empirical distribution* $\mu_N$ associated to $X^1, X^2, \ldots X^N$ is defined by

$$\mu_N(\omega) := \frac{1}{N}\sum_{i=1}^{N}\delta_{X^i}(\omega), \quad \omega \in \Omega. \qquad (2)$$

Notice that (2) defines a map from $\Omega$ to the space $\mathcal{D}$ of probability distributions on $C([0, 1]; \mathbb{R}^n)$. Hence, if $E$ is a subset of $\mathcal{D}$, it makes sense to consider $\mathbb{P}(\omega : \mu_N(\omega) \in E)$. By the ergodic theorem, see e.g. [29, Theorem A.9.3.], the distributions $\mu_N$ converge weakly [1] to $W$ as $N$ tends to infinity. Hence, if $W \notin E$, we must have $\mathbb{P}(\omega : \mu_N(\omega) \in E) \searrow 0$. Large deviation theory provides us with a much finer result: Such a decay is *exponential* and the exponent may be characterized solving a *maximum entropy problem*. Indeed, in our setting, let $E = \mathcal{D}(\rho_0, \rho_1)$, namely distributions on $C([0, 1]; \mathbb{R}^n)$ having marginal densities $\rho_0$ and $\rho_1$ at times $t = 0$ and $t = 1$, respectively. Then, Sanov's theorem, roughly asserts that if the "prior" $W$ does not have

---

[1] Let $\mathcal{V}$ be a metric space and $\mathcal{D}(\mathcal{V})$ be the set of probability measures defined on $\mathcal{B}(\mathcal{V})$, the Borel $\sigma$-field of $\mathcal{V}$. We say that a sequence $\{P_N\}$ of elements of $\mathcal{D}(\mathcal{V})$ converges weakly to $P \in \mathcal{D}(\mathcal{V})$, and write $P_N \Rightarrow P$, if $\int_{\mathcal{V}} f dP_N \to \int_{\mathcal{V}} f dP$ for every bounded, continuous function $f$ on $\mathcal{V}$.

the required marginals, the probability of observing an empirical distribution $\mu_N$ in $\mathcal{D}(\rho_0, \rho_1)$ decays according to

$$\mathbb{P}\left(\frac{1}{N}\sum_{i=1}^{N}\delta_{X^i} \in \mathcal{D}(\rho_0, \rho_1)\right) \sim \exp\left[-N\inf\left\{\mathbb{D}(P\|W); P \in \mathcal{D}(\rho_0, \rho_1)\right\}\right],$$

where

$$\mathbb{D}(P\|W) = \begin{cases} \mathbb{E}_P\left(\log\frac{\mathrm{d}P}{\mathrm{d}W}\right), & \text{if } P \ll W \\ +\infty & \text{otherwise} \end{cases}.$$

is the relative entropy functional or Kullback-Leibler divergence between $P$ and $W$. Thus, the most likely random evolution between two given marginals is the solution of the Schrödinger Bridge Problem:

**Problem 1.**

$$\text{Minimize} \quad \mathbb{D}(P\|W) \quad \text{over} \quad P \in \mathcal{D}(\rho_0, \rho_1). \tag{3}$$

The optimal solution is called the *Schrödinger bridge* between $\rho_0$ and $\rho_1$ over $W$, and its marginal flow $(\rho_t)$ is the *entropic interpolation*.

Let $P \in \mathcal{D}$ be a finite-energy diffusion, namely under $P$ the canonical coordinate process $X_t(\omega) = \omega(t)$ has a (forward) Ito differential

$$dX_t = \beta_t dt + dW_t \tag{4}$$

where $\beta_t$ is adapted to $\{\mathcal{F}_t^-\}$ ($\mathcal{F}_t^-$ is the $\sigma$-algebra of events observable up to time $t$) and

$$\mathbb{E}_P\left[\int_0^1 \|\beta_t\|^2 dt\right] < \infty. \tag{5}$$

Let

$$P_x^y = P\left[\cdot \mid X_0 = x, X_1 = y\right], \quad W_x^y = W\left[\cdot \mid X_0 = x, X_1 = y\right]$$

be the disintegrations of $P$ and $W$ with respect to the initial and final positions. Let also $\pi$ and $\pi^W$ be the joint initial-final time distributions under $P$ and $W$, respectively. Then, we have the following decomposition of the relative entropy [32]

$$\mathbb{D}(P\|W) = E_P\left[\log\frac{dP}{dW}\right] =$$

$$\int\int\left[\log\frac{\pi(x,y)}{\pi^W(x,y)}\right]\pi(x,y)dxdy + \int\int\left(\log\frac{dP_x^y}{dW_x^y}\right)dP_x^y\pi(x,y)dxdy.$$

Both terms are nonnegative. We can make the second zero by choosing $P_x^y = W_x^y$. Thus, the problem reduces to the static one:

**Problem 2.** *Minimize over densities $\pi$ on $\mathbb{R}^n \times \mathbb{R}^n$ the index*

$$\mathbb{D}(\pi \| \pi^W) = \int \int \left[ \log \frac{\pi(x,y)}{\pi^W(x,y)} \right] \pi(x,y) dx dy \tag{6}$$

*subject to the (linear) constraints*

$$\int \pi(x,y) dy = \rho_0(x), \quad \int \pi(x,y) dx = \rho_1(y). \tag{7}$$

If $\pi^*$ solves the above problem, then

$$P^*(\cdot) = \int_{\mathbb{R}^n \times \mathbb{R}^n} W_{xy}(\cdot) \pi^*(x,y) dx dy,$$

solves Problem 1.

Consider now the case when the prior is $W_\gamma$, namely Wiener measure with variance $\gamma$, so that

$$p(0, x, 1, y) = [2\pi\gamma]^{-\frac{n}{2}} \exp \left[ -\frac{|x-y|^2}{2\gamma} \right].$$

Let $\rho_0^W$ denote the initial marginal density under $W_\gamma$. Notice now that the quantity

$$\int \int \left[ \log \rho_0^W(x) \right] \pi(x,y) dx dy = \int \left[ \log \rho_0^W(x) \right] \rho_0(x) dx$$

is independent of $\pi$ satisfying (7) (as long as the support of $\rho_0$ is contained in the support of $\rho_0^W$ the integral is well defined). Using this fact and $\pi^{W_\gamma}(x,y) = \rho_0^W(x) p(0, x; 1, y)$, we now get

$$\begin{aligned}
\mathbb{D}(\pi \| \pi^{\mathbf{W}_\gamma}) &= -\int \int \left[ \log \pi^W(x,y) \right] \pi(x,y) dx dy + \int \int \left[ \log \pi(x,y) \right] \pi(x,y) dx dy \\
&= \int \int \frac{|x-y|^2}{2\gamma} \pi(x,y) dx dy - \mathcal{S}(\pi) + C,
\end{aligned} \tag{8}$$

where $\mathcal{S}$ is the differential entropy and $C$ does not depend on $\pi^2$. Thus, Problem 2 of minimizing $\mathbb{D}(\pi \| \pi^{\mathbf{W}_\gamma})$ over $\Pi(\rho_0, \rho_1)$, namely the "couplings" of $\rho_0$ and $\rho_1{}^3$ is equivalent to

$$\inf_{\pi \in \Pi(\rho_0, \rho_1)} \int \frac{|x-y|^2}{2} \pi(x,y) dx dy + \gamma \int \pi(x,y) \log \pi(x,y) dx dy, \tag{9}$$

---

[2]It follows, in particular, that the initial marginal density of the prior can WLOG always be taken equal to $\rho_0$.

[3]Probability densities on $\mathbb{R}^n \times \mathbb{R}^n$ with marginals $\rho_0$ and $\rho_1$.

namely a regularization of Optimal Mass Transport (OMT) [59] with quadratic cost function obtained by subtracting a term proportional to the entropy.

## 2.3 Derivation of the Schrödinger system

We outline the derivation of the so-called Schrödinger system for the sake of continuity in exposition. Two good surveys on Schrödinger Bridges are [60, 42]. The Lagrangian function for Problem 2 has the form

$$\mathcal{L}(\pi; \lambda, \mu) = \int \int \left[ \log \frac{\pi(x,y)}{\pi^W(x,y)} \right] \pi(x,y) dx dy$$

$$+ \int \lambda(x) \left[ \int \pi(x,y) dy - \rho_0(x) \right] + \int \mu(y) \left[ \int \pi(x,y) - \rho_1(y) \right].$$

Setting the first variation with respect to $\pi$ equal to zero, we get the (sufficient) optimality condition

$$1 + \log \pi^*(x,y) - \log p(0,x,1,y) - \log \rho_0^W(x) + \lambda(x) + \mu(y) = 0,$$

where we have used the expression $\pi^W(x,y) = \rho_0^W(x) p(0,x,1,y)$ with $p$ as in (1). We get

$$\begin{aligned} \frac{\pi^*(x,y)}{p(0,x,1,y)} &= \exp \left[ \log \rho_0^W(x) - 1 - \lambda(x) - \mu(y) \right] \\ &= \exp \left[ \log \rho_0^W(x) - 1 - \lambda(x) \right] \exp \left[ -\mu(y) \right]. \end{aligned}$$

Hence, the ratio $\pi^*(x,y)/p(0,x,1,y)$ factors into a function of $x$ times a function of $y$. Denoting these by $\hat{\varphi}(x)$ and $\varphi(y)$, respectively, we can then write the optimal $\pi^*(\cdot, \cdot)$ in the form

$$\pi^*(x,y) = \hat{\varphi}(x) p(0,x,1,y) \varphi(y), \tag{10}$$

where $\varphi$ and $\hat{\varphi}$ must satisfy

$$\hat{\varphi}(x) \int p(0,x,1,y) \varphi(y) dy = \rho_0(x), \tag{11}$$

$$\varphi(y) \int p(0,x,1,y) \hat{\varphi}(x) dx = \rho_1(y). \tag{12}$$

Let us define $\hat{\varphi}(0,x) = \hat{\varphi}(x), \quad \varphi(1,y) = \varphi(y)$ and

$$\hat{\varphi}(1,y) := \int p(0,x,1,y) \hat{\varphi}(0,x) dx, \quad \varphi(0,x) := \int p(0,x,1,y) \varphi(1,y).$$

Then, (11)-(12) can be replaced by the system

$$\hat{\varphi}(1, y) = \int p(0, x, 1, y)\hat{\varphi}(0, x)dx, \tag{13}$$

$$\varphi(0, x) = \int p(0, x, 1, y)\varphi(1, y)dy, \tag{14}$$

coupled by the boundary conditions

$$\varphi(0, x) \cdot \hat{\varphi}(0, x) = \rho_0(x), \quad \varphi(1, y) \cdot \hat{\varphi}(1, y) = \rho_1(y). \tag{15}$$

Notice that dividing both sides of (10) by $\rho_0(x)$, we get

$$p^*(0, x, 1, y) = \frac{1}{\varphi(0, x)}p(0, x, 1, y)\varphi(1, y), \tag{16}$$

where $\varphi$, in Doob's language, is *space time harmonic* satisfying (14) or, equivalently,

$$\frac{\partial \varphi}{\partial t} + \frac{1}{2}\Delta\varphi = 0. \tag{17}$$

The solution is namely obtained from the prior distribution via a *multiplicative functional transformation* of the prior Markov processes [37]. The question of existence and uniqueness of positive functions $\hat{\varphi}$, $\varphi$ satisfying (13, 14, 15), left open by Schrödinger, is a highly nontrivial one and has been settled in various degrees of generality by Fortet, Beurlin, Jamison and Föllmer [33, 5, 38, 32]. The pair $(\varphi, \hat{\varphi})$ is unique up to multiplication of $\varphi$ by a positive constant $c$ and division of $\hat{\varphi}$ by the same constant. At each time $t$, the marginal $\rho_t$ factorizes as

$$\rho_t(x) = \varphi(t, x) \cdot \hat{\varphi}(t, x). \tag{18}$$

Schrödinger: "Merkwürdige Analogien zur Quantenmechanik, die mir sehr des Hindenkens wert erscheinen"[4] Indeed (18) resembles Born's relation

$$\rho_t(x) = \psi(t, x) \cdot \bar{\psi}(t, x)$$

with $\psi$ and $\bar{\psi}$ satisfying two adjoint equations like $\varphi$ and $\hat{\varphi}$. Moreover, the solution of Problem 1 exhibits the following remarkable *reversibility property*: Swapping the two marginal densities $\rho_0$ and $\rho_1$ the new solution is simply the time reversal of the previous one, cf. the title "On the reversal of natural laws" of [52].

---

[4]Remarkable analogies to quantum mechanics which appear to me very worth of reflection.

### 2.4 "Half bridges"

Consider the following variant of Problem 1 with prior distribution $W_\gamma$:

**Problem 3.**

$$\text{Minimize} \quad \mathbb{D}(P\|W_\gamma) \quad \text{over} \quad P \in \mathcal{D}(\rho_1), \tag{19}$$

namely, we only impose the final marginal. The same argument as before shows that Problem 3 reduces to the following variant of Problem 2:

**Problem 4.** *Minimize over densities $\pi$ on $\mathbb{R}^n \times \mathbb{R}^n$ the index*

$$\mathbb{D}(\pi\|\pi^{W_\gamma}) = \int \int \left[ \log \frac{\pi(x,y)}{\pi^{W_\gamma}(x,y)} \right] \pi(x,y) dx dy \tag{20}$$

*subject to the (linear) constraint*

$$\int \pi(x,y) dx = \rho_1(y). \tag{21}$$

The same variational analysis as in Subsection 2.3, now gives the optimality condition

$$1 + \log \pi^*(x,y) - \log p(0,x,1,y) - \log \rho_0^W(x) + \mu(y) = 0.$$

We then get

$$\frac{\pi^*(x,y)}{p(0,x,1,y)} = \exp\left[\log \rho_0^W(x) - 1 - \mu(y)\right] = \rho_0^W(x) \exp\left[-1 - \mu(y)\right]. \tag{22}$$

Thus, in the previous notation, we can set $\hat{\varphi}(x) = \rho_0^W(x)$ and $\varphi(y) = \exp\left[-1 - \mu(y)\right]$. Let

$$\rho_1^{W_\gamma}(y) = \int \left[2\pi\gamma\right]^{-\frac{n}{2}} \exp\left[-\frac{|x-y|^2}{2\gamma}\right] \rho_0^W(x) dx$$

which replaces (13) with $\hat{\varphi}(0,x) = \rho_0^W(x)$ and $\hat{\varphi}(1,y) = \rho_1^{W_\gamma}(y)$. Then (12) gives immediately

$$\varphi(y) = \frac{\rho_1(y)}{\rho_1^{W_\gamma}(y)}. \tag{23}$$

We now get the form of the optimal initial-final joint distribution of the half-bridge:

$$\begin{aligned}
\pi^*(x,y) &= \rho_0^W(x) p(0,x,1,y) \frac{\rho_1(y)}{\rho_1^{W_\gamma}(y)} \\
&= \rho_0^W(x) \left[2\pi\gamma\right]^{-\frac{n}{2}} \exp\left[-\frac{|x-y|^2}{2\gamma}\right] \frac{\rho_1(y)}{\rho_1^{W_\gamma}(y)} = \pi_\gamma^W(x,y) \frac{\rho_1(y)}{\rho_1^{W_\gamma}(y)}.
\end{aligned}$$

Finally, let

$$\varphi(0, x) := \int (2\pi)^{-\frac{n}{2}} \left[2\pi\gamma\right]^{-\frac{n}{2}} \exp\left[-\frac{|x-y|^2}{2\gamma}\right] \varphi(y) dy. \qquad (24)$$

Then, the initial marginal of the solution is given by

$$\rho_0(x) = \varphi(0, x)\rho_0^W(x).$$

Notice that here there is no delicate question about existence and uniqueness for the Schrödinger system as $\hat{\varphi}$ coincides at all times with the prior one-time marginal. This, in turn, provides the terminal condition for the $\varphi$ function at time $t = 1$ which then only needs to be propagated backward through (24) to provide the full solution. In the special case when $\rho_0^W(x) = \delta(x)$, we have $\rho_t^{W_\gamma}(x) = (2\pi\gamma t)^{-\frac{n}{2}} \exp\left[-\frac{|x|^2}{2t}\right]$ and, in particular, $\rho_1^{W_\gamma}(y) = (2\pi\gamma)^{-\frac{n}{2}} \exp\left[-\frac{|y|^2}{2\gamma}\right]$.

An immediate application of the half-bridge problem is the reconstruction of the past of a system given its current state and a prior model for its evolution. The availability of a prior here is crucial, as without a prior or other regularization such inverse problems are typically ill-posed. Another application concerns deviations from equilibrium. Consider a stochastic system whose dynamics $p(t_1, x_1, t_2, x_2)$ has a statistically steady state $\rho_{eq}(\mathrm{x})$, possibly modulated in time. What is the most likely path that would take us at time $t$ to a state $\rho_1(x)$ away from equilibrium? For example, one may want to anticipate the likely path of strong storms or large waves, so as to be able to forecast them.

## 2.5 Stochastic control and fluid-dynamic formulations

In addition to the formulations above, there exist also dynamic versions of the problem such as the following stochastic control formulation originating with [24, 25, 49]: Problem 1 (when the prior has variance $\gamma$) is equivalent to

**Problem 5.**

$$\text{Minimize}_{u \in \mathcal{U}} \ J(u) = \mathbb{E}\left[\int_0^1 \frac{1}{2\gamma} \|u_t\|^2 dt\right], \qquad (25)$$

$$\text{subject to } dX_t = u_t dt + \sqrt{\gamma} dW_t, \quad X_0 \sim \rho_0(x)dx, \quad X_1 \sim \rho_1(y)dy,$$

where the family $\mathcal{U}$ consists of adapted, finite-energy control functions.

The optimal control is of the feedback type

$$u_t = \gamma \nabla \log \varphi(t, X_t),\tag{26}$$

where $(\varphi, \hat\varphi)$ solve the Schrödinger system (13, 14, 15). These formulations are particularly relevant in applications where the prior distribution on paths is not simply the Wiener measure, but is associated to the uncontrolled ("free") evolution of a dynamical system, see e.g [16, 17, 14] and in image morphing/interpolation [13, Subsection 5.3]. In the case of the half bridge, (26) still holds with $\varphi$ satisfying

$$\frac{\partial\varphi}{\partial t} + \frac{\gamma}{2}\Delta\varphi = 0, \quad \varphi(1, \cdot) = \frac{\rho_1(\cdot)}{\rho_1^{\mathbf{W}^\gamma}(\cdot)}.$$

Problem 5 leads immediately to the following fluid dynamic problem:

**Problem 6.**

$$\inf_{(\rho, b)} \int_{\mathbb{R}^n} \int_0^1 \frac{1}{2}\|b(x, t)\|^2 \rho(t, x) dt dx,\tag{27a}$$

$$\frac{\partial\rho}{\partial t} + \nabla\cdot(b\rho) - \frac{\gamma}{2}\Delta\rho = 0,\tag{27b}$$

$$\rho(0, x) = \rho_0(x), \quad \rho(1, y) = \rho_1(y).\tag{27c}$$

*where $b(\cdot, \cdot)$ varies over continuous functions on $\mathbb{R}^n \times [0, 1]$.*

This problem is not equivalent to Problems 1, 2 and 5 in that it only reproduces the optimal *entropic interpolating flow* $\{\rho_t; 0 \le t \le 1\}$. Information about correlations at different times and smoothness of the trajectories is here lost. As $\gamma \searrow 0$, the solution to this problem converges to the solution of the Benamou-Brenier Optimal Mass Transport problem [3, 44, 45, 46, 42, 41]:

$$\inf_{(\rho, v)} \int_{\mathbb{R}^n} \int_0^1 \frac{1}{2}\|v(x, t)\|^2 \rho(t, x) dt dx,\tag{28a}$$

$$\frac{\partial\rho}{\partial t} + \nabla\cdot(v\rho) = 0,\tag{28b}$$

$$\rho(0, x) = \rho_0(x), \quad \rho(1, y) = \rho_1(y).\tag{28c}$$

Let $(\rho, b)$ be optimal for Problem 6 and define the *current velocity field* [47]

$$
\begin{aligned}
v(x, t) &= b(x, t) - \frac{\gamma}{2}\nabla\log\rho_t(x)\\
&= \gamma\nabla\log\varphi(t, x) - \frac{\gamma}{2}\nabla\log\rho_t(x) = \frac{\gamma}{2}\nabla\log\frac{\varphi(t, x)}{\hat\varphi(t, x)},
\end{aligned}\tag{29}
$$

where we have used (26) and (18). Assume that $v$ guarantees existence and uniqueness of the initial value problem on $[0, 1]$ for any deterministic initial condition and consider

$$\dot{X}(t) = v(X(t), t), \quad X(0) \sim \rho_0 dx. \tag{30}$$

Then the probability density $\rho_t(x)$ of $X(t)$ satisfies (weakly) the continuity equation

$$\frac{\partial \rho}{\partial t} + \nabla \cdot (v\rho) = 0$$

as well as (27b) with the same initial condition and therefore coincides with $\rho(x, t)$. This suggests that an alternative fluid-dynamic problem characterizing the entropic interpolation flow $\{\rho_t; 0 \leq t \leq 1\}$ may be possible. Indeed, such time-symmetric problem was derived in [15]:

**Problem 7.**

$$\inf_{(\rho, v)} \int_{\mathbb{R}^n} \int_0^1 \left[ \frac{1}{2} \|v(x, t)\|^2 + \frac{\gamma}{8} \|\nabla \log \rho\|^2 \right] \rho(t, x) dt dx, \tag{31a}$$

$$\frac{\partial \rho}{\partial t} + \nabla \cdot (v\rho) = 0, \tag{31b}$$

$$\rho(0, x) = \rho_0(x), \quad \rho(1, y) = \rho_1(y). \tag{31c}$$

The two criteria differ by $(\gamma/8)\mathcal{I}(\rho)$ where the Fisher information functional $\mathcal{I}$ is given by

$$\mathcal{I}(\rho) = \int \|\nabla \log \rho_t\|^2 \rho_t(x) dx$$

while the Fokker-Planck equation (27b) has been replaced by the continuity equation (31b). Both Problems 6 and 7 can be thought of as regularizations of the Benamou-Brenier problem (28) and as dynamic counterparts of (9). Also notice that, precisely as in Problem (28), the optimal current velocity (29) in Problem 7 is of the gradient type.

Finally, consider the family of diffeomorphisms $\{T_t; 0 \leq t \leq 1\}$ satisfying

$$\frac{dT_t}{dt}(x) = v(T_t(x), t), \quad T_0 = I, \tag{32}$$

where $v$ is defined by (29). Then, in analogy to the *displacement interpolation* of Optimal Mass Transport, we have the following relation for the entropic interpolation flow

$$\rho_t(x) dx = T_t \# \rho_0(x) dx, \tag{33}$$

namely $\rho_t(x)dx$ is the *push-forward* of the measure $\rho_0(x)dx$ under the map $T_t$. In particular, the map $T_\gamma = T_1$ pushes $\rho_0(x)dx$ onto $\rho_1(x)dx$ and represents therefore the entropic counterpart of the map solving the original Monge problem. It may be called the *Monge-Schrödinger map*.

## 3  Fortet's iterative algorithm

The oldest proof of existence and uniqueness for the Schrödinger system (13, 14, 15), due to Fortet [33], is *algorithmic* in nature, establishing convergence of successive approximations. More explicitly, let $g(x, y)$ be a nonnegative, continuous function bounded from above. Suppose $g(x, y) > 0$ except possibly for a zero measure set for each fixed value of $x$ or of $y$. Suppose that $\rho_0(x)$ and $\rho_1(y)$ are continuous, nonnegative functions such that

$$\int \rho_0(x)dx = \int \rho_1(y)dy.$$

Suppose, moreover, that the integral

$$\int \frac{\rho_1(y)}{\int g(z, y)\rho_0(z)dz} dy$$

is finite. Then, [33, Theorem 1], the system

$$\phi(x) \int g(x, y)\psi(y)dy = \rho_0(x), \tag{34}$$

$$\psi(y) \int g(x, y)\phi(x)dx = \rho_1(y) \tag{35}$$

admits a solution $(\phi(x), \psi(y))$ with $\phi \geq 0$ continuous and $\psi \geq 0$ measurable. Moreover, $\phi(x) = 0$ only where $\rho_0(x) = 0$ and $\psi(y) = 0$ only where $\rho_1(y) = 0$.

The result is proven by setting up a complex approximation scheme to show that equation

$$h(x) = \Omega(h) = \int g(x, y) \frac{\rho_1(y)dy}{\int g(z, y)\frac{\rho_0(z)}{h(z)}dz}. \tag{36}$$

has a positive solution. Notice that

$$g(x, y) = p(0, x, 1, y) = [2\pi\gamma]^{-\frac{n}{2}} \exp\left[-\frac{|x - y|^2}{2\gamma}\right].$$

satisfies all assumptions of Fortet's theorem. Uniqueness, in the sense described in Subsection 2.3, namely uniqueness of rays, is much easier to establish. In the recent paper [30], the bulk of Fortet's paper has been rewritten filling in all the missing steps and providing explanations for the rationale behind the various articulations of his approach.

Independently, at about the same time and in the discrete setting, an *iterative proportional fitting* (IPF) procedure, was proposed in the statistical literature on contingency tables [26]. Convergence for the IPF algorithm was first established (in a special case) by Richard Sinkhorn in 1964 [54]. The iterates were shortly afterwards shown to converge to a "minimum discrimination information" [36, 31], namely to a minimum entropy distance. This line of research, usually called *Sinkhorn algorithms*, continues to this date, see e.g. [23, 1, 57].

It is apparent that an iterative scheme can be designed based on (36) which, in the previous notation, reads

$$\Omega(\varphi(0, x)) = \int p(0, x, 1, y) \frac{\rho_1(y) dy}{\int p(0, z, 1, y) \frac{\rho_0(z)}{\varphi(0, z)} dz}. \tag{37}$$

This was accomplished in [13], showing convergence of the iterates in a suitable projective metric, but only for the case when both marginals have compact support. We outline the approach below. Let $\mathcal{S}$ be a real Banach space and let $\mathcal{K} \subset \mathcal{S}$ be a closed cone with nonempty interior int$\mathcal{K}$ and such that $\mathcal{K} + \mathcal{K} \subseteq \mathcal{K}$, $\mathcal{K} \cap -\mathcal{K} = \{0\}$ as well as $\lambda\mathcal{K} \subseteq \mathcal{K}$ for all $\lambda \geq 0$. Define the partial order

$$x \preceq y \Leftrightarrow y - x \in \mathcal{K}, \quad x < y \Leftrightarrow y - x \in \text{int}\mathcal{K}$$

and for $x, y \in \mathcal{K}_0 := \mathcal{K}\backslash\{0\}$, define

$$\begin{aligned} M(x, y) &:= \inf\{\lambda \mid x \preceq \lambda y\} \\ m(x, y) &:= \sup\{\lambda \mid \lambda y \preceq x\}. \end{aligned}$$

Then, the Hilbert metric is defined on $\mathcal{K}_0 = \mathcal{K}\backslash\{0\}$ by

$$d_H(x, y) := \log\left(\frac{M(x, y)}{m(x, y)}\right).$$

Strictly speaking, it is a *projective* metric since it is invariant under scaling by positive constants, i.e., $d_H(x, y) = d_H(\lambda x, \mu y)$ for any $\lambda > 0, \mu > 0$ and $x, y \in \text{int}\mathcal{K}$. Thus, it is actually a distance between *rays*. Important examples in finite dimension are provided by the positive orthant of $\mathbb{R}^n$ and by the cone of positive semidefinite matrices of dimension $n \times n$. One can then take advantage of some remarkable contractivity theorems on cones established for suitable positive linear and nonlinear maps most noticeably by Garret Birkhoff and P. Bushell [6, 11]. We mention that the celebrated Perron-Frobenius theorem may be viewed as a corollary of these results [7].

Let $S_i \subset \mathbb{R}^n$, $i = 0, 1$ be the compact support of $\rho_i$, $i = 0, 1$. In order to study the Schrödinger system (13, 14, 15) we consider the maps

$$\mathcal{E}: \ \varphi(1, y) \mapsto \varphi(0, x) = \int_{S_1} p(0, x, 1, y) \varphi(1, y) dy \qquad (38a)$$

$$\mathcal{E}^\dagger: \ \hat\varphi(0, x) \mapsto \hat\varphi(1, y) = \int_{S_0} p(0, x, 1, y) \hat\varphi(0, x) dx \qquad (38b)$$

$$\mathcal{D}_0: \ \varphi(0, x) \mapsto \hat\varphi(0, x) = \rho_0(x) / \varphi(0, x) \qquad (38c)$$

$$\mathcal{D}_1: \ \hat\varphi(1, y) \mapsto \varphi(1, y) = \rho_1(y) / \hat\varphi(1, y), \qquad (38d)$$

on appropriate domains. The map $\Omega$ in (37), can then be written as the composition

$$\Omega = \mathcal{E} \circ \mathcal{D}_1 \circ \mathcal{E}^\dagger \circ \mathcal{D}_0.$$

While maps $\mathcal{D}_0$ and $\mathcal{D}_1$ can be seen to be isometries in the Hilbert metric, the two linear maps are always non-expanding and, under suitable assumptions, strictly contractive. Things are complicated by the fact that certain cones in infinite dimensions such as the nonnegative integrable functions, have nonempty interior and the contractivity theorems do not apply. It is, however, possible to represent the map $\Omega$ as an alternative composition of other four maps all acting on a cone with nonempty interior so that contractivity of $\Omega$ can be established, see [13, Section 3] for the details. Once convergence of the rays is established in the projective metric, it suffices to use the fact that both $\varphi(0, x) \cdot \hat\varphi(x, 0)$ and $\varphi(1, y) \cdot \hat\varphi(1, y)$ must integrate to one to show convergence of the functions.

Setting up an iterative scheme based on (37) when only samples of the two marginals are available is obviously much more challenging: This is the main topic of this paper which we shall pursue starting from the next section. This will also provide an approach to data-driven Optimal Mass Transport alternative to [58] since, as observed at the end of Subsection 2.2, the Schrödinger Bridge Problem may be viewed as a regularization of OMT.

# 4 Numerical methodology

This section develops a sample-based numerical methodology for the solution of the Schrödinger bridge problem. This is the case, ubiquitous in applications, where the distributions $\rho_0$ and $\rho_1$ are only known through the finite sample sets $\{x_i\}$ and $\{y_j\}$ of cardinality $m$ and $n$ respectively.

One could propose a scheme whereby one first estimates $\rho_0$ and $\rho_1$ from the samples provided, and then solves the regular Schrödinger bridge problem between these two estimates. Yet there are a number of reasons why a procedure based directly on the sample sets is preferable:

1. Density estimation adds an extra computational layer to the algorithm, and hence a source of additional potential approximation errors.

2. In high-dimensional settings, density estimation is inherently difficult and requires larger data sets than are customarily available.

3. Even with estimations for $\rho_0$ and $\rho_1$ known in closed form, the solution to the Schrödinger bridge problem requires the calculation of integrals that in most cases cannot be performed in closed form. Then, rather than introducing grids, whose size grows exponentially with the dimension of the space, and whose local mesh-sizes are difficult to adapt to the local distributions, it is better to resort to samples, which provide a naturally adapted discretization of the continuous problem, and which yield a Monte-Carlo error that scales only mildly with the dimensionality of the space.

For conciseness, we shall denote $p(y|x)$ the prior transition density $p(0, x, 1, y)$, and write $\hat{\varphi}_0(x)$ and $\varphi_1(y)$ instead of $\hat{\varphi}(0, x)$ and $\varphi(1, y)$, respectively. Then (10) reads

$$\pi^*(x, y) = \hat{\varphi}_0(x) \ p(y|x) \ \varphi_1(y).$$

The entropic interpolation between $\rho_0$ and $\rho_1$ is given by $\rho_t(z) = \varphi_t(z)\hat{\varphi}_t(z)$, where

$$\varphi_t(z) = \int p(t, z, 1, y)\varphi_1(y)dy \qquad \hat{\varphi}_t(z) = \int p(0, x, t, z)\hat{\varphi}_0(x)dx. \quad (39)$$

In particular, one needs to solve the system

$$\rho_0(x) = \varphi_0(x)\hat{\varphi}_0(x), \quad \rho_1(y) = \varphi_1(y)\hat{\varphi}_1(y),$$

with

$$\varphi_0(x) = \int p(y|x)\varphi_1(y)dy \qquad \hat{\varphi}_1(y) = \int p(y|x)\hat{\varphi}_0(x)dx.$$

To begin, we need to reformulate the problem so that it involves the distributions only through their available samples.

## 4.1 The half-bridge problem through maximal likelihood

We develop first an algorithm for the half-bridge problem. Even though this is much simpler than the full bridge, it includes some of its main ingredients. The data-driven version of equation (23) for the half-bridge problem is: given $\hat{\varphi}_1(y) \geq 0$ and a set of $n$ samples $\{y_j\}$ of $\rho_1(y)$, find $\varphi_1(y)$ such that $\hat{\varphi}_1(y)\varphi_1(y) = \rho_1(y)$, an equality that we interpret in the maximum likelihood sense:

$$\varphi_1 = \underset{\varphi_1(y) \geq 0}{\arg\max} \sum_j \log\left(\hat{\varphi}_1(y_j)\varphi_1(y_j)\right), \quad \text{subject to } \int \left(\hat{\varphi}_1(y)\varphi_1(y)\right) dy = 1.$$

Consider first a situation with infinitely many samples $y_j$, or equivalently with $\rho_1(y)$ known. Then the problem becomes

$$\varphi_1 = \underset{\varphi_1(y) \geq 0}{\arg\max} \int \log\left(\varphi_1(y)\right) \rho_1(y) \, dy, \quad \int \left(\hat{\varphi}_1(y)\varphi_1(y)\right) dy = 1.$$

We can satisfy the positivity constraint automatically by proposing an exponential form for $\varphi_1$:

$$\varphi_1(y) = e^{g(y)},$$

which yields

$$\max_g \int g(y)\rho_1(y)dx \quad \text{s.t.} \quad \int \hat{\varphi}_1(y)e^{g(y)}dy = 1,$$

or, introducing a Lagrange multiplier $\lambda$ for the constraint,

$$\max_g \min_\lambda L(g, \lambda) = \int g(y)\rho_1(y) \, dx - \lambda \left( \int \hat{\varphi}_1(y)e^{g(y)}dy - 1 \right).$$

Maximizing over $g$ first yields

$$\frac{\delta L}{\delta g} = \rho_1(y) - \lambda\hat{\varphi}_1(y)e^{g(y)} = 0 \rightarrow g(y) = \log\left( \frac{\rho_1(y)}{\lambda\hat{\varphi}_1(y)} \right).$$

Then the minimization over $\lambda$ becomes

$$\min_\lambda \left[ -\log(\lambda) + \lambda \right] \Rightarrow \lambda = 1.$$

Hence the value of the optimal $\lambda$ is known explicitly, and the estimation problem becomes:

$$\max_g L(g) = \int g(y)\rho_1(y)\ dx - \int \hat{\varphi}_1(y)e^{g(y)}dy + 1. \tag{40}$$

Notice that the solution to (40) is

$$g(y) = \log\left(\frac{\rho_1(y)}{\hat{\varphi}_1(y)}\right) \Rightarrow \varphi_1(y) = \frac{\rho_1(y)}{\hat{\varphi}_1(y)},$$

the exact answer to the problem. Yet in the true problem $\rho_1(y)$ is only known through samples $\{y_j\}$, so the first integral in (40) must be replaced by its empirical counterpart:

$$\int g(y)\rho_1(y)\ dx \to \frac{1}{n}\sum_j g(y_j).$$

Then, introducing a rough estimate $\tilde{\rho}_1$ of $\rho_1$ that one can sample, such as a Gaussian, and drawing $\tilde{n}$ samples $\tilde{y}_k$ from it, we can replace the second integral above by its Monte Carlo simulation:

$$\int \hat{\varphi}_1(y)e^{g(y)}dy = \int \frac{\hat{\varphi}_1(y)e^{g(y)}}{\tilde{\rho}_1(y)}\tilde{\rho}_1(y)dy \to \frac{1}{\tilde{n}}\sum_k \frac{\hat{\varphi}_1(\tilde{y}_k)e^{g(\tilde{y}_k)}}{\tilde{\rho}_1(\tilde{y}_k)}.$$

(Notice that for $\tilde{\rho}_1 = \rho_1$ and $g$ the true maximizer, this is an estimation with zero variance.)

Finally, proposing a parameterization of the unknown $g(y)$, such as

$$g(y) = \sum_l \beta_l F_l(y),$$

where the $F_l$ are functions externally provided, we end up with the following algorithm for estimating $\varphi_1(y)$:

$$\varphi_1(y) = e^{\sum_l \beta_l F_l(y)},$$

where $\beta$ solves

$$\beta = \arg\max L = \sum_l \left(\frac{1}{n}\sum_j F_l(y_j)\right)\beta_l - \frac{1}{\tilde{n}}\sum_k \frac{\hat{\varphi}_1(\tilde{y}_k)e^{\sum_l \beta_l F_l(\tilde{y}_k)}}{\tilde{\rho}_1(\tilde{y}_k)},$$

a convex optimization problem, since the

$$\frac{\partial^2 L}{\partial \beta_i \beta_j} = -\frac{1}{\tilde{n}} \sum_k \frac{\hat{\varphi}_1(\tilde{y}_k) e^{\sum_l \beta_l F_l(\tilde{y}_k)}}{\tilde{\rho}_1(\tilde{y}_k)} F_i(\tilde{y}_k) F_j(\tilde{y}_k)$$

form a negative definite matrix.

More generally, we could have adopted a form for $\varphi_1(y) = \Phi(y, \beta)$ different from the exponential, while still guaranteeing positivity, such as

$$\Phi(y, \beta) = g(y, \beta)^2,$$

where $g(y, \beta)$ is any family of real functions with parameters $\beta$. Then the problem above would have become

$$\beta = \arg \max L = \frac{1}{n} \sum_j \log \left( \Phi(y_j, \beta) - \frac{1}{\tilde{n}} \sum_k \frac{\hat{\varphi}_1(\tilde{y}_k) \Phi(\tilde{y}_k, \beta)}{\tilde{\rho}_1(\tilde{y}_k)} \right).$$

## 4.2 The full bridge problem

Since the solution of the Schrödinger problem is given in (10) by

$$\pi^*(x, y) = \hat{\varphi}_0(x) p(y|x) \varphi_1(y),$$

it is natural to parameterize in closed form only the functions $\hat{\varphi}_0(x)$ and $\varphi_1(x)$. As in the half-bridge problem, we guarantee the positivity of these two functions directly through their parameterization $\hat{\varphi}_0(x, \hat{\beta}), \varphi_1(x, \beta)$, for instance writing them as the exponential or square of some other real functions.

If $\hat{\varphi}_1$ were given, we would find the coefficients $\beta$ defining $\varphi_1$ by solving an optimization problem entirely analogous to the half-bridge problem before:

$$\beta = \arg \max L_1 = \frac{1}{n} \sum_j \log \left( \varphi_1(y_j, \beta) \right) - \int \hat{\varphi}_1(y) \varphi_1(y, \beta) dy.$$

However, at every step in the algorithm, only $\hat{\varphi}_0(x)$ is available in closed form; in order to find $\hat{\varphi}_1(y)$ we need to propagate the former through

$$\hat{\varphi}_1(y) = \int p(y|x) \hat{\varphi}_0(x, \hat{\beta}) dx.$$

Then

$$\int \hat{\varphi}_1(y) \varphi_1(y, \beta) dy = \int \left[ \int p(y|x) \hat{\varphi}_0(x, \hat{\beta}) dx \right] \varphi_1(y, \beta) dy$$

$$= \int \left[ \int p(y|x) \varphi_1(y, \beta) dy \right] \hat{\varphi}_0(x, \hat{\beta}) dx.$$

Since the inner integral equals $\varphi_0(x)$, and $\varphi_0(x)\hat{\varphi}_0(x) = \rho_0(x)$, we can multiply and divide by a sampleable estimator $\tilde{\rho}_0$ of $\rho_0$ with $\tilde{m}$ samples $\{\tilde{x}_i\}$, and write

$$\int \hat{\varphi}_1(y)\varphi_1(y,\beta)dy \approx \frac{1}{\tilde{m}} \sum_i \left[ \int p(y|\tilde{x}_i)\varphi_1(y,\beta)dy \right] \frac{\hat{\varphi}_0(\tilde{x}_i,\hat{\beta})}{\tilde{\rho}_0(\tilde{x}_i)},$$

an estimation with zero variance at the exact solution if $\tilde{\rho}_0 = \rho_0$. Since the $\tilde{x}_i$ are fixed throughout the algorithm, we can at little expense extract, for each $i$, $\hat{n}$ samples $\hat{y}_i^j$ from the prior $p(y|\tilde{x}_i)$, and write the final estimator

$$\int \hat{\varphi}_1(y)\varphi_1(y,\beta)dy \approx \frac{1}{\tilde{m}\hat{n}} \sum_{i,j} \varphi_1(\hat{y}_i^j,\beta) \frac{\hat{\varphi}_0(\tilde{x}_i,\hat{\beta})}{\tilde{\rho}_0(\tilde{x}_i)},$$

so the problem for $\beta$ becomes

$$\beta = \arg\max \frac{1}{n} \sum_j \log\left(\varphi_1(y_j,\beta)\right) - \frac{1}{\tilde{m}\hat{n}} \sum_{i,j} \varphi_1(\hat{y}_i^j,\beta) \frac{\hat{\varphi}_0(\tilde{x}_i,\hat{\beta})}{\tilde{\rho}_0(\tilde{x}_i)}. \qquad (41)$$

For the parameters $\hat{\beta}$, we have

$$\hat{\beta} = \arg\max L_0 = \frac{1}{m} \sum_i \log\left(\hat{\varphi}_0(x_i,\hat{\beta})\right) - \int \hat{\varphi}_0(x,\hat{\beta})\varphi_0(x)dx,$$

where

$$\varphi_0(x) = \int p(y|x)\varphi_1(y,\beta)dy.$$

Then

$$\begin{aligned} \int \hat{\varphi}_0(x,\hat{\beta})\varphi_0(x)dx &= \int \left[ \int p(y|x)\varphi_1(y,\beta)dy \right] \hat{\varphi}_0(x,\hat{\beta})dx \\ &\approx \frac{1}{\tilde{m}\hat{n}} \sum_{i,j} \varphi_1(\hat{y}_i^j,\beta) \frac{\hat{\varphi}_0(\tilde{x}_i,\hat{\beta})}{\tilde{\rho}_0(\tilde{x}_i)}. \end{aligned}$$

(The fact that this is exactly the same estimation than for the integral $\int \hat{\varphi}_1(y)\varphi_1(y,\beta)dy$ should not be entirely surprising, as both equal one and involve the same parameters.) Finally,

$$\hat{\beta} = \arg\max \frac{1}{m} \sum_i \log\left(\varphi_0(x_i,\hat{\beta})\right) - \frac{1}{\tilde{m}\hat{n}} \sum_{i,j} \varphi_1(\hat{y}_i^j,\beta) \frac{\hat{\varphi}_0(\tilde{x}_i,\hat{\beta})}{\tilde{\rho}_0(\tilde{x}_i)}. \qquad (42)$$

### 4.3 The algorithm

Summarizing the results above, we have developed the following algorithm:

1. **Data:** We are provided with $m$ samples $\{x_i\}$ of $\rho_0(x)$, $n$ samples $\{y_j\}$ of $\rho_1(y)$, and a prior conditional probability density $p(y|x)$. The latter needs not be known in closed form, but one should be able to sample it for any value of $x$ (if the opposite is true, i.e. we know $p(y|x)$ in closed form but cannot sample it, an alternative algorithm presented below should be applied.)

2. **Goal:** To find the most likely joint distribution $\pi(x, y)$ under the prior $p(y|x)$ consistent with the two marginals, and the corresponding posterior $p^*(y|x)$. When $p(y|x)$ is the end result of the prior $p(t_1, x, t_2, y)$ for a time dependent process, we also seek the more detailed posterior $p^*(t_1, x, t_2, y)$ for this process, as well as the intermediate distributions $\rho_t(z)$ for $t \in [0, 1]$.

3. **Preliminary work:** Based on the samples $\{x_i\}$, we need to produce a first estimate $\tilde{\rho}_0$ of $\rho_0(x)$ and $\tilde{m}$ independent samples $\{\tilde{x}_i\}$ drawn from it. More specifically, we will need these $\tilde{m}$ samples and the values $\tilde{\rho}_0(\tilde{x}_i)$ of $\tilde{\rho}_0$ on them. For instance, one can use the Gaussian kernel density estimator

$$\tilde{\rho}_0(x) = \frac{1}{m} \sum_i G(x - x_i),$$

where $G$ is an isotropic Gaussian with suitable bandwidth. In building this estimate, we can use, in addition to the samples $\{x_i\}$, any additional prior information that we may have on $\rho_0(x)$. For instance, its support may be known to be contained within some set $\Omega$, typically not to include unrealistic negative values of some components of $x$. One simple way to address this particular case is to multiply the unconstrained estimator $\tilde{\rho}_0$ by the characteristic function of $\Omega$, reject any sample outside of $\Omega$, and normalize the resulting distribution through division by the factor

$$\frac{\tilde{m} + m_r}{\tilde{m}},$$

where $m_r$ is the total number of rejections that occurred.

For each sample $\tilde{x}_i$, we need to produce $\hat{n}$ samples $\hat{y}_i^j$ drawn independently from $p(y|\tilde{x}_i)$. For instance, if $p$ is the result of a diffusive process

between $t = 0$ and $t = 1$, with drift $u(x, t)$ and diffusivity $\nu(x, t)$, we would simulate the stochastic process

$$dx = u(x, t)dt + \nu(x, t)dW, \quad x(0) = x_i, \quad y_i^j = x(1).$$

If $p(y|x)$ is known in closed form but is not easily sampled, one can propose another conditional probability $q(y|x)$ not very far from $p$ but sampleable, and produce weighted samples $y_i^j$ from $q(y|\tilde{x}_i)$, with weights

$$w_i^j = \frac{p(y_i^j|\tilde{x}_i)}{q(y_i^j|\tilde{x}_i)},$$

to be included as extra factors under the second sum in problems (41) and (42).

4. **Model selection and initialization:** We need to propose a parametric family of non-negative real functions $\Phi(z, \beta)$. Examples are

$$\Phi(z, \beta) = e^{\sum_k \beta_k F_k(z)} \quad \text{and} \quad \Phi(z, \beta) = \left( \sum_k \beta_k F_k(z) \right)^2, \quad (43)$$

where the $F_k$ are a given set of functions (monomials, Legendre functions, sines and cosines, splines, etc.) In high-dimensions, we may want to use instead a low-rank tensor factorization as in [35, 55] The final estimated joint density will adopt the form

$$\pi(x, y) = \Phi(x, \hat{\beta}) \ p(y|x) \ \Phi(y, \beta),$$

and the estimated posterior conditional probability will be

$$P^*(y|x) = \frac{p(y|x) \ \Phi(y, \beta)}{\int p(z|x) \ \Phi(z, \beta)dz},$$

where the integral in the denominator can be estimated for each desired value of $x$ by simulating $p(z|x)$. In the notation above,

$$\hat{\varphi}_0(x) = \Phi(x, \hat{\beta}) \quad \text{and} \quad \varphi_1(y) = \Phi(y, \beta).$$

We initialize the algorithm with an initial guess for $\beta$, such as the $\beta$ that yields the default $\varphi_1(y) = 1$ (i.e. $\beta = 0$ when using the first of the parametrizations in (43). This is typically easier than starting with a guess for $\hat{\beta}$ approximating the corresponding default $\hat{\varphi}_0(x) = \rho_0(x)$). When using the quadratic parametrization in (43), we start with a choice of $\hat{\beta}$ that, depending on the chosen basis functions $F_l$, yields to the biggest effective support of $\phi(x)$.

5. **Main loop:** We alternate between the updates (42) for $\hat{\beta}$ and (41) for $\beta$ iteratively until a convergence criterion is met. Some choices for the family $\Phi(z, \beta)$, such as

$$\Phi(z, \beta) = e^{\sum_k \beta_k F_k(z)}$$

yield automatically convex optimization problems for $\hat{\beta}$ and $\beta$.

# 5 Numerical examples

This section illustrates the proposed methodology on two examples relevant in applications: the interpolation of probability distributions, and a variation on importance sampling in the context of Monte Carlo estimates of integrals.

## 5.1 Interpolation between two Gaussian mixtures

Figure 1 displays the two marginal distributions of a two dimensional numerical example, where $\rho_0$ and $\rho_1$ are Gaussian mixtures given by

$$\rho_0 = \frac{1}{3} \sum \left[ \mathcal{N}(\mu_1, \Sigma_1) + \mathcal{N}(\mu_2, \Sigma_2) + \mathcal{N}(\mu_3, \Sigma_3) \right]$$

$$\rho_1 = \frac{1}{3} \left[ \mathcal{N}(\mu_4, \Sigma_4) + \mathcal{N}(\mu_5, \Sigma_5) + \mathcal{N}(\mu_6, \Sigma_6) \right]$$

with parameters

$$\mu_1 = \begin{bmatrix} -2 \\ 1.5 \end{bmatrix}, \Sigma_1 = \begin{bmatrix} 0.2 & 0.1 \\ 0.1 & 0.4 \end{bmatrix}, \mu_2 = \begin{bmatrix} 0.2 \\ 1.2 \end{bmatrix}, \Sigma_2 = \begin{bmatrix} 0.6 & -0.4 \\ -0.4 & 0.6 \end{bmatrix},$$

$$\mu_3 = \begin{bmatrix} 0.5 \\ -1 \end{bmatrix}, \Sigma_3 = \begin{bmatrix} 0.5 & 0.4 \\ 0.4 & 0.7 \end{bmatrix}, \mu_4 = \begin{bmatrix} -1.8 \\ 1.1 \end{bmatrix}, \Sigma_4 = \begin{bmatrix} 0.3 & 0.1 \\ 0.1 & 0.3 \end{bmatrix},$$

$$\mu_5 = \begin{bmatrix} -0.2 \\ 1.2 \end{bmatrix}, \Sigma_5 = \begin{bmatrix} 0.5 & -0.3 \\ -0.3 & 0.8 \end{bmatrix} \mu_6 = \begin{bmatrix} -0.5 \\ 0.9 \end{bmatrix}, \Sigma_6 = \begin{bmatrix} 0.6 & 0.2 \\ 0.2 & 0.6 \end{bmatrix}. \quad (44)$$

Figure 2 displays the interpolation between $\rho_0$ and $\rho_1$ obtained by computing $\rho_t(z) = \varphi_t(z)\hat{\varphi}_t(z)$ for each time $t \in [0, 1]$ at the data points $z(t)$ obtained by integrating the equation (29, 30) with $\gamma = 2$. In this example, both $\varphi$ and $\hat{\varphi}$ were represented as the square of linear combinations of the first 10 Hermite functions.

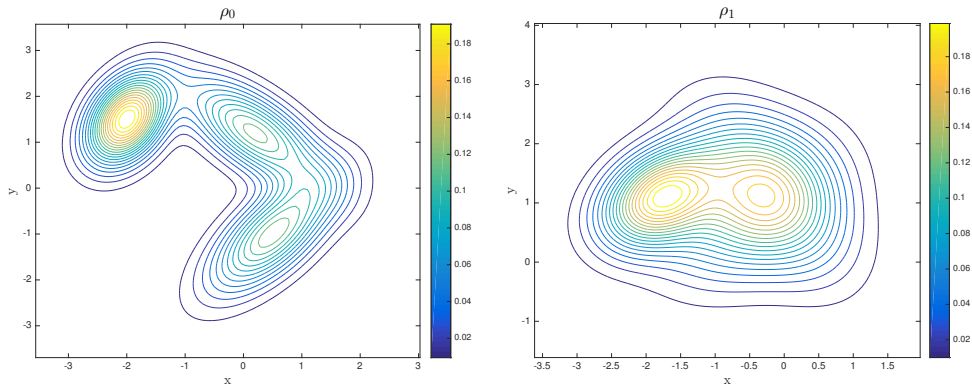

Figure 1: Initial and final probability density distribution from which points $x_i$ and $y_j$ where sampled respectively. This is the only input used by the algorithm.

## 5.2 A variation on importance sampling

The methodology of this article turns out to be particularly well suited to improve Monte Carlo estimates of the quantity

$$I = \int f(y)\rho_1(y)dy, \tag{45}$$

when $\rho_1(y)$ is only known through $n$ sample points drawn from it. It is know that ordinary Monte Carlo estimates suffer of a slow convergence rate as a function of $n$. Moreover, when the support of $f$ is localized in regions where the value of $\rho_1$ is small, we may have very few points where $f$ is substantially different from zero. If $\rho_1$ where known in closed form, we could remedy these problems though importance sampling, whereby we would rewrite (45) in the form

$$I = \int f(y)\rho_1(y)dy = \int \frac{f(y)\rho_1(y)}{\mu(y)}\mu(y)dy,$$

where $\mu(y)$ is a distribution easy to sample and such that $f\rho_1/\mu$ has small variance, and then estimate $I$ via Monte Carlo:

$$I \approx \frac{1}{n}\sum_{i=1}^{n}\frac{f(z_i)\rho_1(z_i)}{\mu(z_i)},$$

where the $z_i$ are samples drawn independently from $\mu$. Yet this procedure requires the capacity to evaluate $\rho_1$ at the given points. We are considering instead the frequently occurring situation where $\rho_1$ is only known through $n$ samples $\{y_j\}$.

In this case, we propose to use the sample points $\{y_j\}$ to solve the Schrödinger bridge problem between $\rho_1(y)$ and a distribution $\rho_0(x)$ of our choice. This allows us to map arbitrary points in $y$-space to $x$-space. In particular, we can chose points $\tilde{y}_j$ that resolve $f$ well, and use them to estimate the integral $I$ through the following steps:

1. Select points $\tilde{y}_j$ spanning the support of $f$. These points do not need to be sampled from a probability distribution, the only requirement is that $f(y)$ be well characterized by its values on the $\tilde{y}_j$.

2. Compute $\varphi_1$ and $\hat{\varphi}_0$ solving the Schrödinger bridge between $\rho_1$ and a standard normal distribution $\rho_0 = \mathcal{N}(0,1)$, through the procedure described in section 4.3. Then integrate the equation of motion (29, 30) by evaluating $\varphi_t$ and $\hat{\varphi}_t$ via Monte Carlo estimates of the integrals in (39) obtained by sampling the prior $p(x,t,y,s)$. Let $T$ be the one-to-one map determined by the solution of (29, 30).

3. Integrate back in time (29, 30) in order to map the points $\tilde{y}_j$ to the points $x_i = T^{-1}(\tilde{y}_i)$ in the domain of $\rho_0$ .

4. Perform a rough Gaussian mixture density estimation $\nu(x)$ of the distribution underlying the points $x_i$ and sample $N$ new points $\tilde{x}_i$ from it.

5. Map the points $\tilde{x}_i$ back to the support of $\rho_1$ to obtain new points $\hat{y}_j = T(\tilde{x}_j)$.

The integral in (45) is then estimated through

$$\int f(y)\rho_1(y)dy = \int f(y(x))\frac{\rho_0(x)}{\nu(x)}\nu(x)dx \approx \frac{1}{N}\sum_{i=1}^{N} f(T(\tilde{x}_j))\frac{\rho_0(\tilde{x}_j)}{\nu(\tilde{x}_j)}. \quad (46)$$

Hence in a sense we have transferred importance sampling from $y$ to the auxiliary $x$-space.

An alternative to the procedure proposed above would first estimate $\rho_1$ using the samples $y_j$, and then estimate $I$ though regular importance sampling. However, estimating $\rho_1$ in high dimensions is challenging, especially when few data points $y_j$ are available. Moreover, the estimation of $\rho_1$ would be particularly poor in those areas where $f$ is large, since $\rho_1$ is small

$$I_{\mathrm{R}} = 0.09894$$
$$I_{\mathrm{MC}} = 0.10512 \pm 0.28865$$
$$I_{\mathrm{S}} = 0.09710 \pm 0.13284$$

Table 1: $I_{\mathrm{R}}$ indicates the reference value for $I$, $I_{\mathrm{MC}}$ is the Monte Carlo estimates of $I$ and $I_{\mathrm{S}}$ is the estimate of $I$ obtained with the procedure described above.

there, hence the number of local samples available will be small. By contrast, in our procedure, the points $y_j$ are used to estimate the (inverse) map from $\rho_1$ to a standard normal distribution instead. Because of the relation $\rho_0(x) = \rho_1(y(x))|\det J(x)|$ the map from $\rho_0$ to $\rho_1$ is in general smoother than $\rho_1$. Therefore, it is in general more robust to parametrize the map rather than the density.

In the numerical experiment in Figure 3, we chose $\rho_1$ to be the equal weight mixture of the three Gaussian: $\mathcal{N}(-1.4, 0.8^2)$, $\mathcal{N}(2.2, 0.4^2)$, $\mathcal{N}(0.2, 0.1^2)$, and $f(y)$ a mixture of the two Gaussian $\mathcal{N}(-0.8, 0.02^2)$, $\mathcal{N}(1, 0.03^2)$, again with equal weights. We compute the reference value $I_{\mathrm{R}}$ for the integral $I = \int f(y)\rho_1(y)dy$ using a uniform grid of step size $h = 10^{-4}$ and compare this value with plain MC estimates of $I$ obtained with 1000 points sampled from $\rho_1$ and with our procedure. As it can be seen from Table 1, the procedure described above gives a better estimates in terms of both the error with respect the reference value and the uncertainty associated with the estimate.

## 6    Conclusions

In this article, we have posed the sample-based Schrödinger bridge problem and developed a methodology for its numerical solution. Characterizing the initial and final distributions of the bridge in terms of samples is well-suited for applications and also natural from a theoretical perspective, since what is a large-deviation problem for a large but finite set of particles becomes a true impossibility as the number of particles grows unboundedly. One must distinguish though between the sample-based formulation, where $\{x_i\}$ and $\{y_j\}$ are regarded as samples of underlying distributions $\rho_0$ and $\rho_1$, from the discrete Schrödinger problem, where the latter are replaced by the empirical distributions $\frac{1}{m}\sum_{i=1}^{m}\delta(x - x_i)$ and $\frac{1}{n}\sum_{j=1}^{n}\delta(y - y_j)$. This article studies the former, finding the joint distribution $\pi^*(x, y)$ for all values of $(x, y)$, not just the sample points, and characterizing the intermediate distributions

$\rho_t(z)$ also for all $z$.

The methodology of this article mimics the iterative scheme developed for the classical bridge problem, but replacing some of its key ingredients by data analogues. Thus the boundary conditions at $t = 0$ and $t = 1$ are reinterpreted in a maximum likelihood sense, thus giving rise to optimization problems, and the integrals defining the propagation of the two factors of $\rho_t$ are estimated via importance sampling.

The data-based Schrödinger problem has a broad scope of applicability. Potential applications include the estimation of atmospheric winds and oceanic currents from tracers, the solution of inverse diffusive problems, the reconstruction of the intermediate evolution of species between well-documented stages, and many more. Since this article is concerned with the development of a general methodology, we have not dwelled into any application in particular, but just illustrated the procedure with two relatively simple examples.

## Acknowledgments

M.P. would like to thank the Courant Institute of Mathematical Sciences of the New York University for the hospitality during the time this research was carried out. Tabak's work was partially supported by NSF grant DMS-1715753 and ONR grant N00014-15-1-2355.

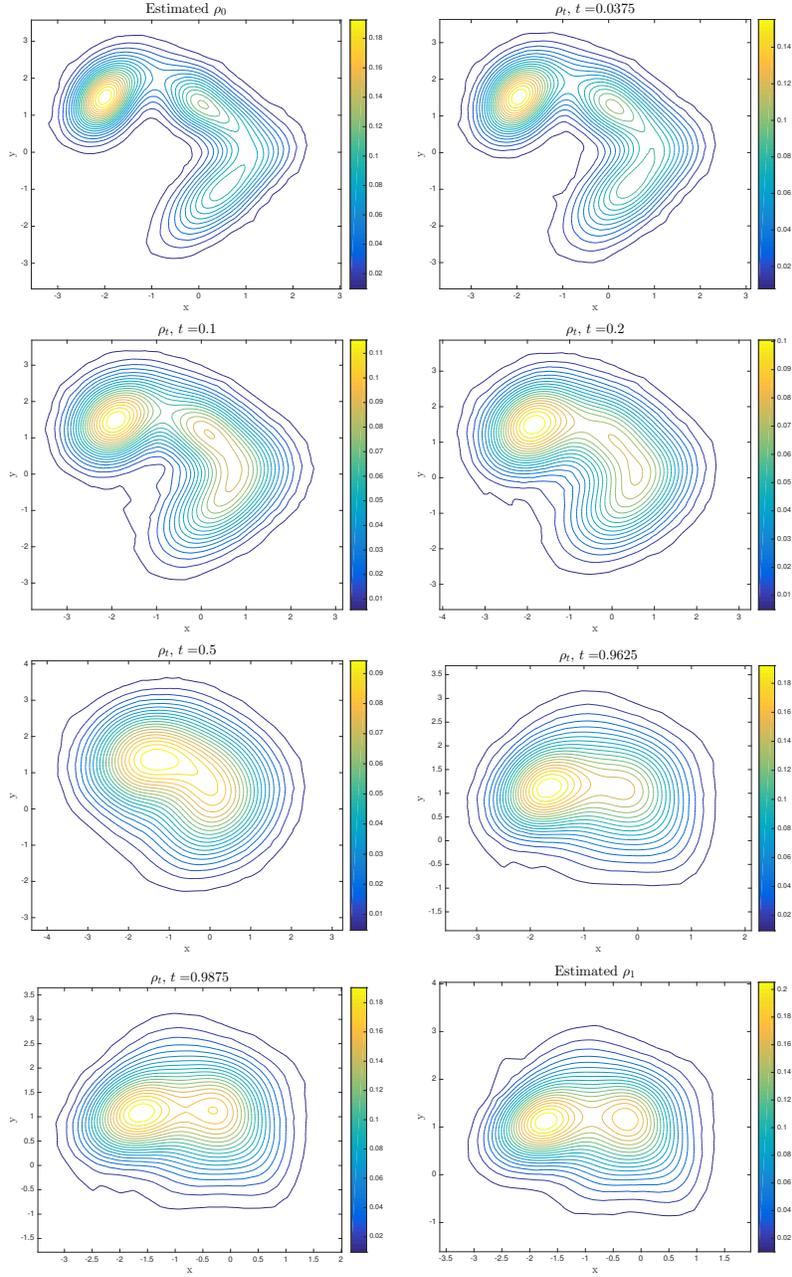

Figure 2: Interpolation between $\rho_0$ and $\rho_1$. Each image is obtained by interpolating $\rho_t(z)$ on the points $z(t)$ representing the solution of (29, 30). Both $\varphi$ and $\hat{\varphi}$ were represented as the square of linear combinations of the first 10 Hermite functions.

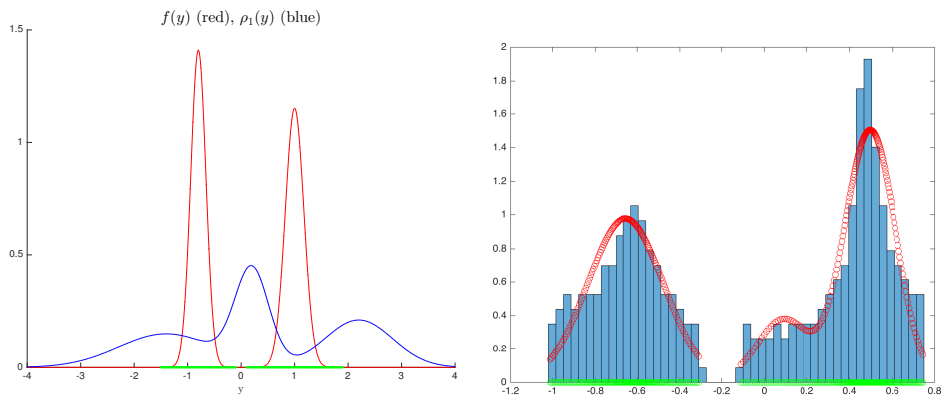

Figure 3: Left panel: The density $\rho_1(y)$ is plotted in blue while the function $f(y)$ is plotted in red. Notice that the support of $f(y)$ is substantially different from zero where the two local minima of $\rho_1$ are placed. The green points $y_j$, appearing one the $x$ axis are points on a regular grid that were selected based on the value of $f$ being bigger than a certain threshold. Right panel: The green points on the left panel have been mapped into the green points $x_j = T^{-1}(y_i)$ on the right panel, in blue there if the histogram of the points $x_i$ and in red its (Gaussian mixture) estimate of it.